\documentclass[12pt]{amsart}
\usepackage{amsthm}
\usepackage{epsfig}

\newtheorem{thm}{Theorem}
\newtheorem{lemma}[thm]{Lemma}

\theoremstyle{definition}

\def\interval#1#2{I_{#1}(#2)}
\def\Itop#1{\interval{top}{#1}}
\def\Ibot#1{\interval{bot}{#1}}
\def\citeBMR{\cite{bogart:trapezoid}}

\def\fig#1#2#3#4{
\begin{figure}[ht]
\begin{center}
\epsfig{file=#1,width=#3}
\caption{#2}
\label{#4}
\end{center}
\end{figure}
}

\def\fign#1#2#3{
\begin{figure}[ht]
\begin{center}
\epsfig{file=#1,width=#3}
\caption{#2}
\end{center}
\end{figure}
}

\begin{document}

\title{Trapezoid Order Classification}
\author{Stephen P. Ryan}
\thanks{Supported by ONR contract N0014-94-1-0950}
\thanks{Research at MSRI is supported in part by NSF grant DMS-9022140}

\address{\hskip-\parindent 
Department of Mathematics, 6188 Bradley Hall, Dartmouth College, 
Hanover, NH 03755}
\email{stephen.ryan@dartmouth.edu}
\date{30 September 1996}

\begin{abstract}In this paper we show the nonequivalence of
combinations of several natural geometric restrictions on trapezoid
representations of trapezoid orders.  Each of the properties unit
parallelogram, unit trapezoid and proper parallelogram, unit trapezoid
and parallelogram, unit trapezoid, proper parallelogram, proper
trapezoid and parallelogram, proper trapezoid, parallelogram, and
trapezoid is shown to be distinct from each of the others.
Additionally, interval orders are shown to be both unit trapezoid and
proper parallelogram orders.
\end{abstract}

\maketitle

\section*{Introduction} A recent paper by Bogart, M\"ohring, and
Ryan \citeBMR\ established that the class of proper
trapezoid orders properly contains the class of unit trapezoid orders.
In this paper we shall make use of the techniques introduced in
that paper, using them to produce a more extensive classification of
trapezoid orders.  

Throughout the following, a {\bf trapezoid order} is an order with a
representation by trapezoids with bases on two parallel lines, called
the baselines, such that $x \prec y$ if and only if the trapezoid
associated to $x$ lies to the left of the trapezoid associated to $y$.
Equivalently, a trapezoid order is an order with interval dimension at
most two.  A {\bf proper} representation is a representation in which no
trapezoid is properly contained in any other.  A {\bf unit} representation
is one in which each trapezoid has the same area.  A {\bf parallelogram}
representation is one in which the upper and lower intervals of each
trapezoid have the same length; i.e. each trapezoid is actually a
parallelogram.  A {\bf proper parallelogram} representation is a
parallelogram representation in which no parallelogram is properly
contained in any other.  A {\bf unit parallelogram} representation is a
parallelogram representation in which each parallelogram has the same
area.  For more detailed descriptions of these definitions, see \citeBMR.

A {\bf rectangle} order is a trapezoid order in which every trapezoid
is actually a rectangle; rectangle orders are actually interval
orders, since then all of the information is contained in the induced
intervals on just one baseline.

We will classify orders according to the existence of any of the
representations described above.  e.g. an order with a proper
parallelogram representation will be referred to as a proper
parallelogram order, even if the current representation we are looking
at is not a proper parallelogram representation.  

The results developed also apply to trapezoid graphs, using the same
techniques to deal with comparability invariance as those used in
\citeBMR.

Trapezoid graphs are the intersection graphs of a set of trapezoids as
defined above.  They were first discussed in \cite{dagan88:trapezoid};
the natural ordering associated with the complement of such graphs,
along with many other properties equivalent to being a trapezoid
graph, is discussed in \cite{langley:tolerance}.

One such property is having interval dimension at most 2.  In this
way, the property of being a trapezoid order (graph) is a
generalization of the property of being an interval order (graph).
The notions of proper and unit trapezoid orders are natural
generalizations of the corresponding concepts for interval orders.
Fred Roberts proved that proper interval is equivalent to unit
interval in \cite{roberts69:indifference}.  \citeBMR\ showed that this
equivalence did not generalize to trapezoid orders.

For convenience, we will restate the jaw lemma from
\citeBMR\ here.

\fig{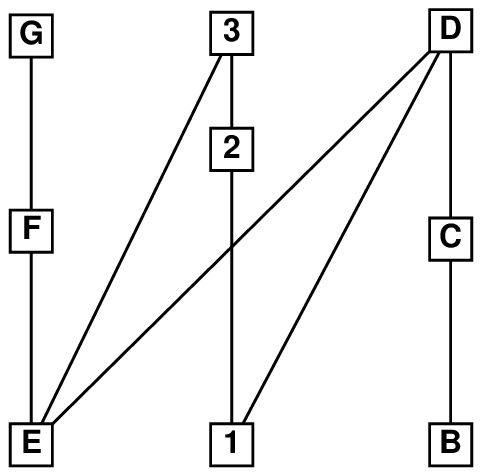}{The jaw order}{3in}{fig:jaworder}

\begin{lemma}[Jaw Lemma]  
\label{lem:jaw}
The order shown in Figure \ref{fig:jaworder} has a trapezoid
representation, and hence is a trapezoid order.  Further, in every
trapezoid representation, we must have endpoints in the relations
$$r(B)<l(C)\leq r(1)<l(2)\leq r(E)<l(D)\leq r(2)<l(3)\leq r(F)<l(G)$$
and
$$R(E)< L(2)\leq R(2)< L(D)$$or in the relations
$$R(B)<L(C)\leq R(1)<L(2)\leq R(E)<L(D)\leq R(2)<L(3)\leq R(F)<L(G)$$
and
$$r(E)< l(2)\leq r(2)< l(D)$$
\footnote{Order drawings produced using Graphlet, a toolkit for implementing
graph editors and graph drawing algorithms.  Graphlet is produced by
the Design, Analysis, Implementation, and Evaluation of Graph Drawing
Algorithms project of the German Science Foundation (DFG) and is
available from the University of Passau on the WWW at
http://www.uni-passau.de/Graphlet/}
\end{lemma}

\fig{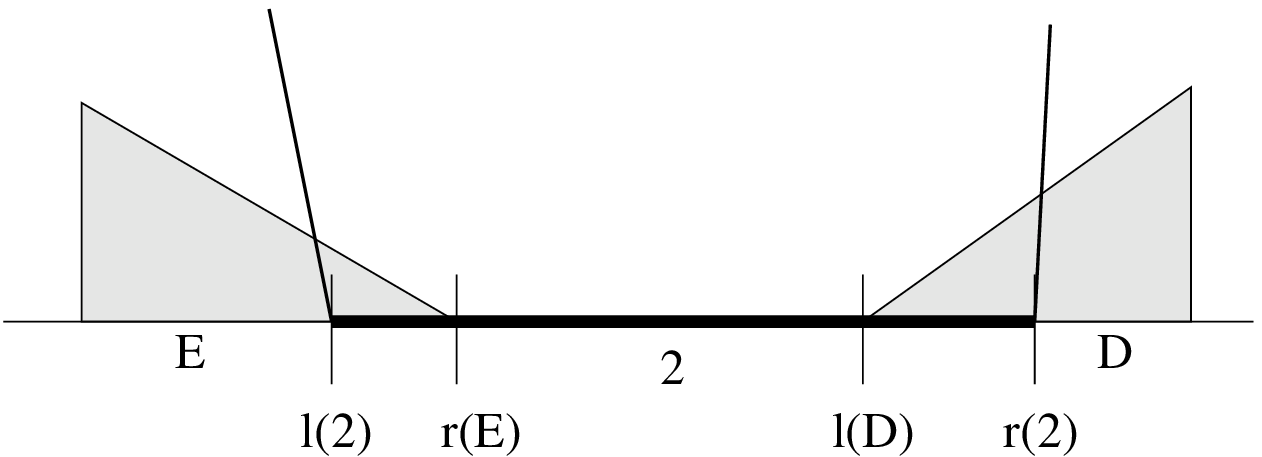}{The ``jaws'' of the jaw lemma}{\hsize}{fig:jaws}

These classes of posets are themselves a poset under the inclusion
relation, and so we may use an order diagram or a Venn diagram to
represent the various inclusions.  A Venn diagram is given in Figure
\ref{fig:diagram}.  Some of the inclusions are obvious; since
parallelograms are special cases of trapezoids, any of the classes of
parallelogram orders are contained in the corresponding class of
trapezoid orders.  Additionally, any class of unit orders is contained
in the corresponding proper orders.  This immediately establishes all
of the inclusions illustrated in Figure \ref{fig:diagram}.  The fact
that all of these inclusions are proper inclusions is the major result
of this paper, generally done by exhibiting an example of an order in
each class not contained in any class below it.

\section*{Preliminaries}

\fig{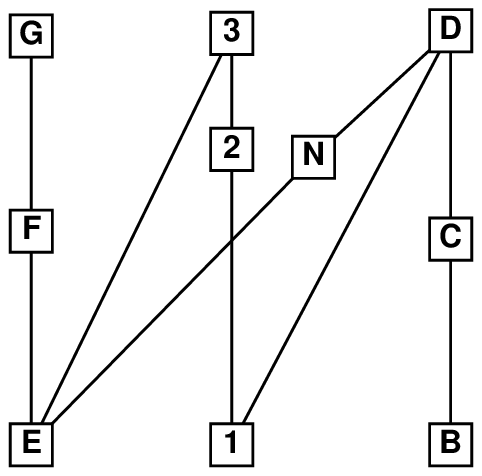}{A unit trapezoid order that is not a
parallelogram order} {3in}{fig:nonparallel}
\begin{thm}\label{thm:nonparallel}
The order in Figure \ref{fig:nonparallel} is a unit trapezoid order,
but has no parallelogram representation.
\end{thm}

\begin{proof} The order in Figure \ref{fig:nonparallel} is the order
in Figure \ref{fig:jaworder} with one additional point ($N$) added.
The jaw lemma applies to the restriction isomorphic to the jaw order,
and therefore we are forced into the trapezoid representation of
Figure \ref{fig:finaljaw}.
\fig{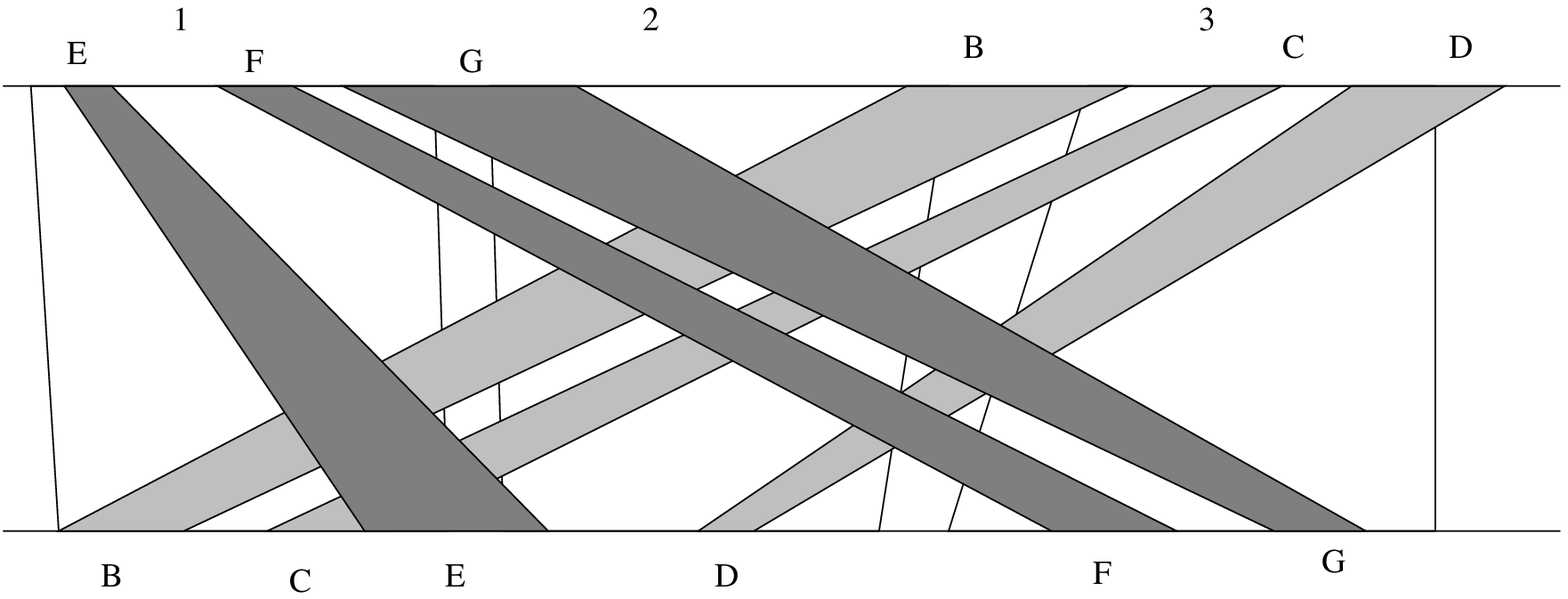}{``The'' trapezoid representation of the jaw
order}{\hsize}{fig:finaljaw}
\fign{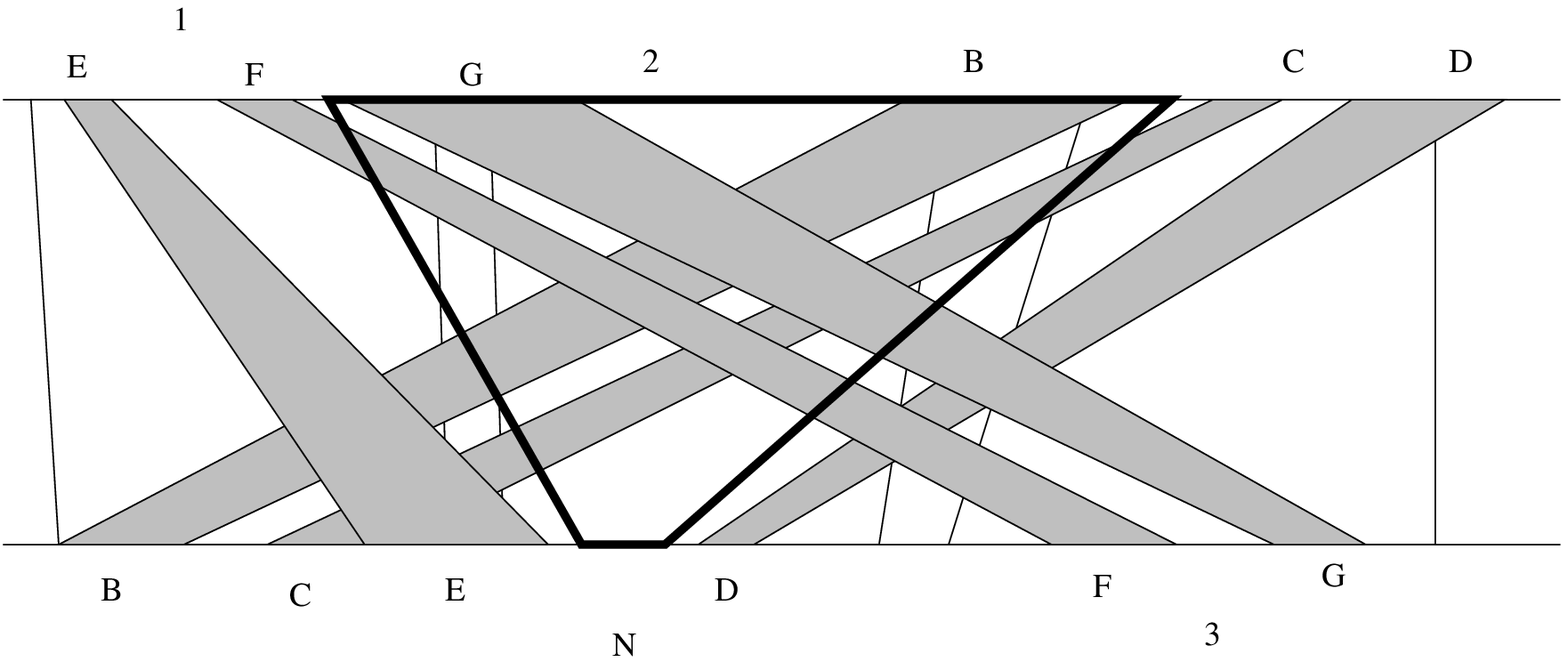}{The trapezoid outlined by the bold lines is
a generic trapezoid representing $N$}{\hsize}

$N$ is between the two teeth of the jaw lemma, and therefore on the
lower baseline, the interval for $N$ is strictly contained in the
interval for 2.  However, $N$ is incomparable to both 1 and 3.
Therefore, the upper left endpoint of $N$ lies to the left of the
upper right endpoint of 1; i.e. $L(N)\leq R(1)$; similarly, $L(3)\leq R(N)$.
However, because $1<2<3$ is a chain, $L(N)<L(2)\leq R(2)<R(N)$, i.e.  the
upper interval for $N$ strictly contains the upper interval for 2.  If
2 is represented by a parallelogram, then $N$ cannot be, and hence,
the order has no parallelogram representation.
\end{proof}

\fig{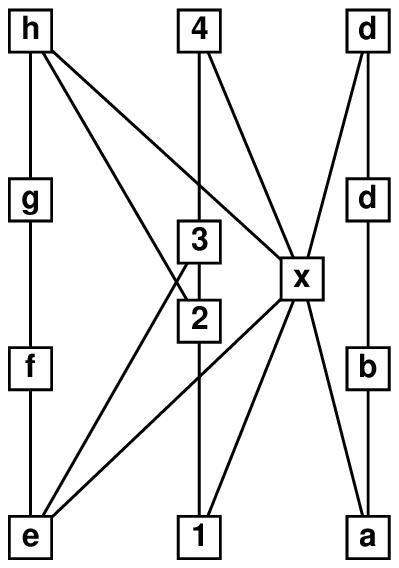}{This order is proper parallelogram, and unit trapezoid,
but not unit parallelogram.}{.4\linewidth}{fig:pput}
\begin{thm}\label{thm:pput}
The order in Figure \ref{fig:pput} is a proper parallelogram order,
and a unit trapezoid order, but not a unit parallelogram order.
\end{thm}

\begin{proof}
To prove that the order in Figure \ref{fig:pput} is a proper
parallelogram order, it is sufficient to exhibit a proper
parallelogram representation.  This is provided in Figure
\ref{fig:pputtrap}.  
\fig{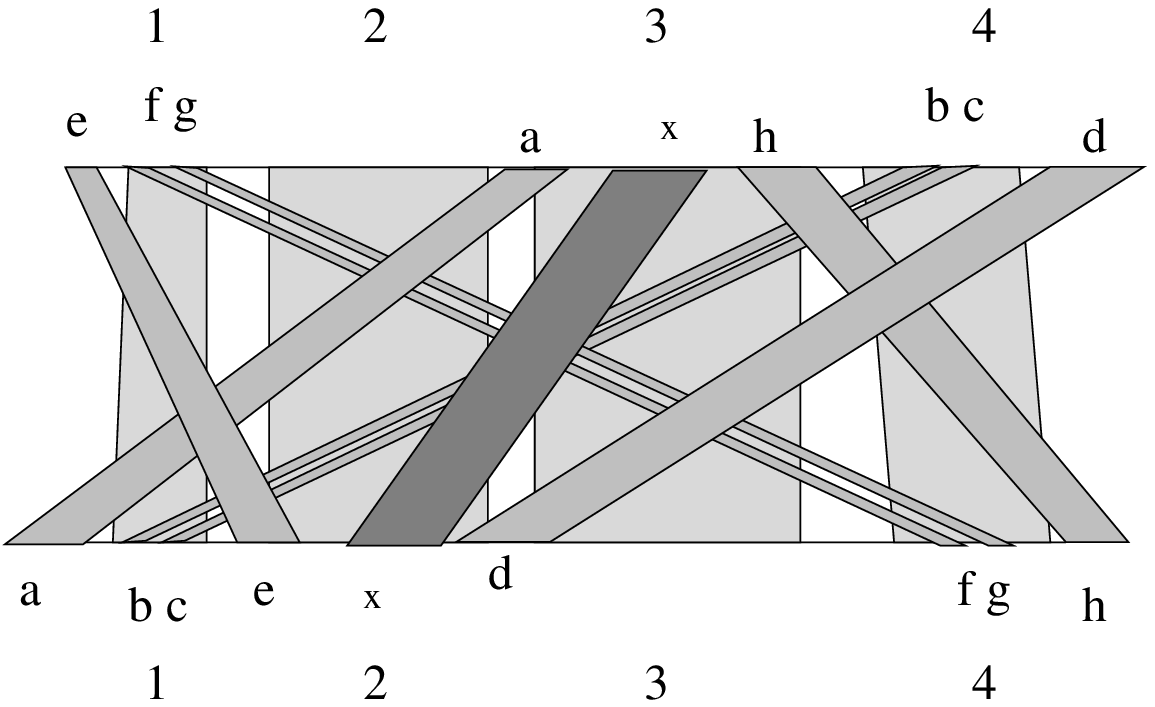}{A parallelogram representation of Figure
\ref{fig:pput}}{\linewidth,height=2.5in}{fig:pputtrap}

To prove that it is a unit trapezoid order, it is sufficient to
exhibit a unit trapezoid representation.  What is shown in Figure
\ref{fig:pputunit} is not a unit trapezoid representation; $b$, $c$,
$f$, and $g$ are all clearly smaller in area than the other
trapezoids.  However, this figure can be modified by moving the
intervals for these trapezoids out and stretching them to the desired
widths.  Note that the bases for $b,c,f$ and $g$ need not lie within
$T_1$ or $T_4$, and can be moved out to the extremities of the diagram.

\fig{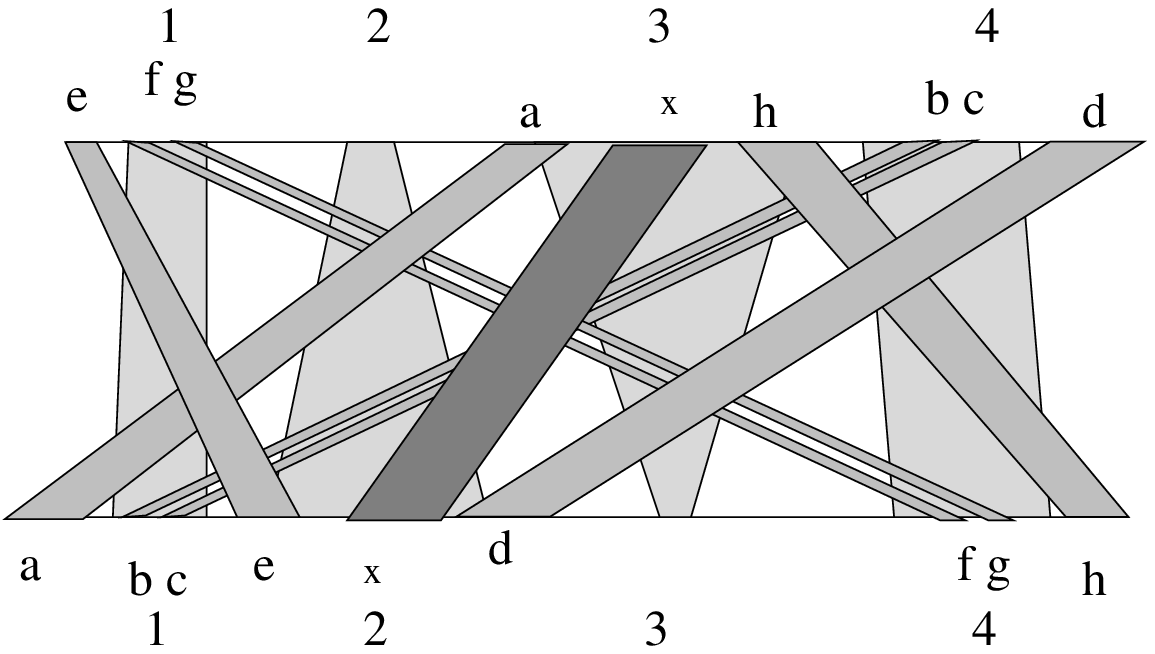}{An approximation of a unit representation of Figure
\ref{fig:pput}}{\linewidth}{fig:pputunit}

Most of this order looks the same as the order in Figure 13 of
\citeBMR.  In fact, the restriction of this order to
$\{a,b,c,d,e,f,g,h,1,2,3,4\}$ is the same as the restriction of Figure
13 of \citeBMR\ to the same elements.  Therefore we get
for free the framework of Figure 16 in \citeBMR, which is
reproduced here in Figure \ref{fig:framework}.
\fig{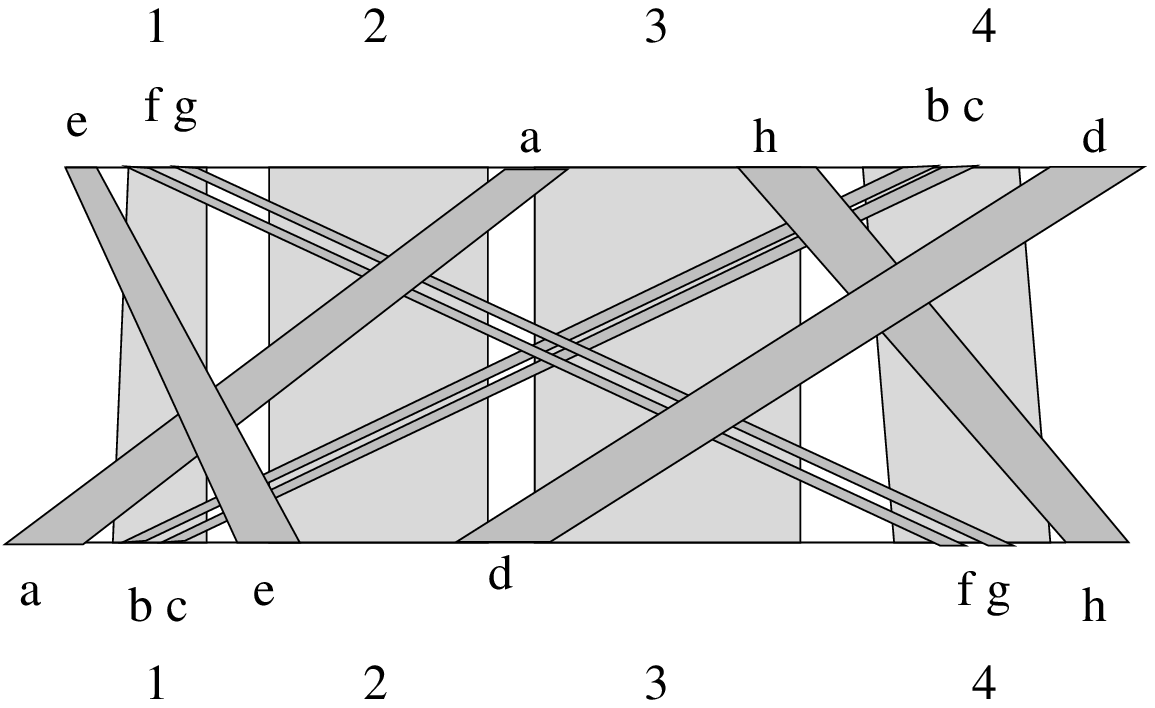}{The framework of Figure 16 in 
\citeBMR.}{\linewidth}{fig:framework}

To show that no parallelogram representation can be a unit
representation, it is sufficient to observe that the upper interval
$\Itop{x}$ must be contained in the upper interval $\Itop{3}$, because
$x$ is contained by the jaws of $a$ and $h$.  If any
representation is to be a parallelogram representation, then the area
of $T_x$ must be strictly less than the area of $T_3$, and therefore
at least one of the two areas cannot be unit.
\end{proof}

\fig{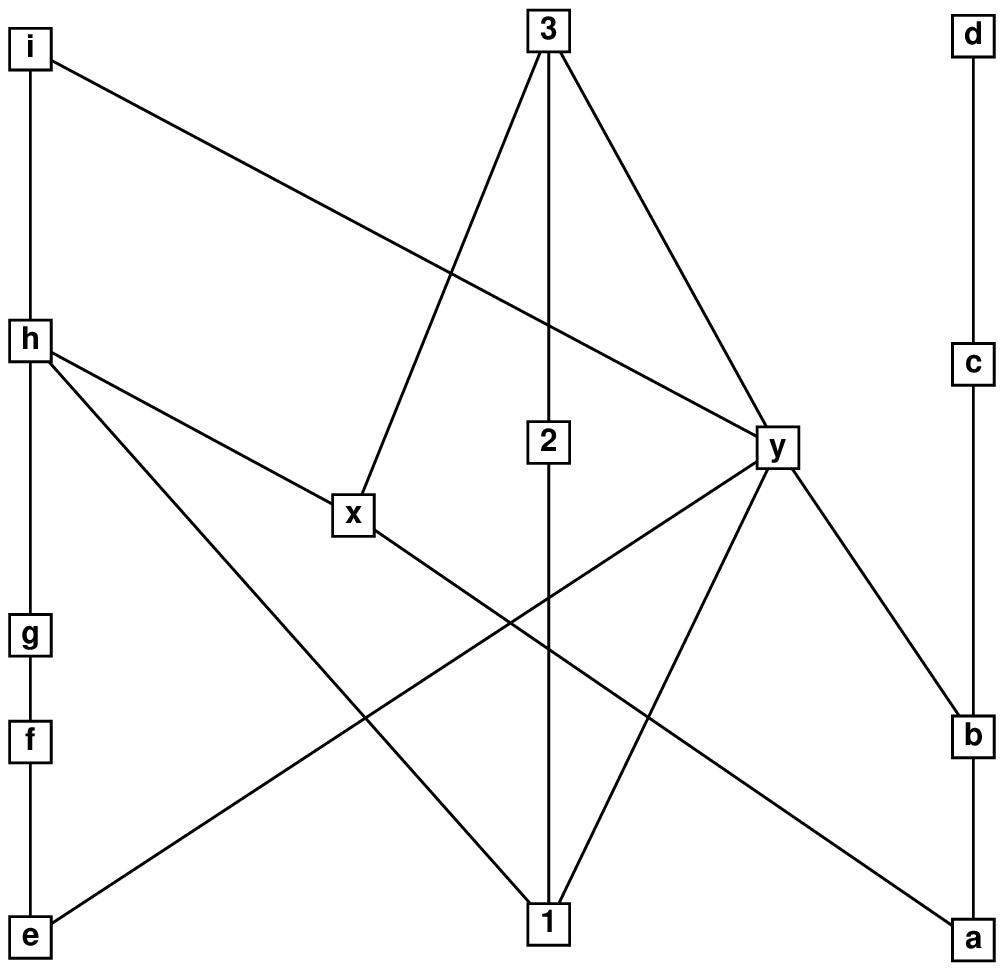}{This is an example of an order that has a unit
trapezoid representation and a parallelogram representation but no
proper parallelogram representation.}{.75\hsize}{fig:uglymess}

\begin{thm}\label{thm:uglymess}
The order in Figure \ref{fig:uglymess} is a unit trapezoid order and a
parallelogram order, but not a proper parallelogram order.
\end{thm}

\begin{proof}This example is probably the most difficult to understand
of all the examples utilized in the classification done in this paper,
due to the subtleness of the property we are trying to establish, and
the fact that we are using two and a half jaw structures.

The statement of this theorem asserts that there is a proper (actually
unit) trapezoid representation, and a parallelogram representation,
but that the two properties cannot be true simultaneously.  To prove
parallelogram, and to prove unit trapezoid, we must exhibit
representations of each.  These are given in Figures \ref{fig:npput},
\ref{fig:npput-p} and \ref{fig:npput-u}.

\fig{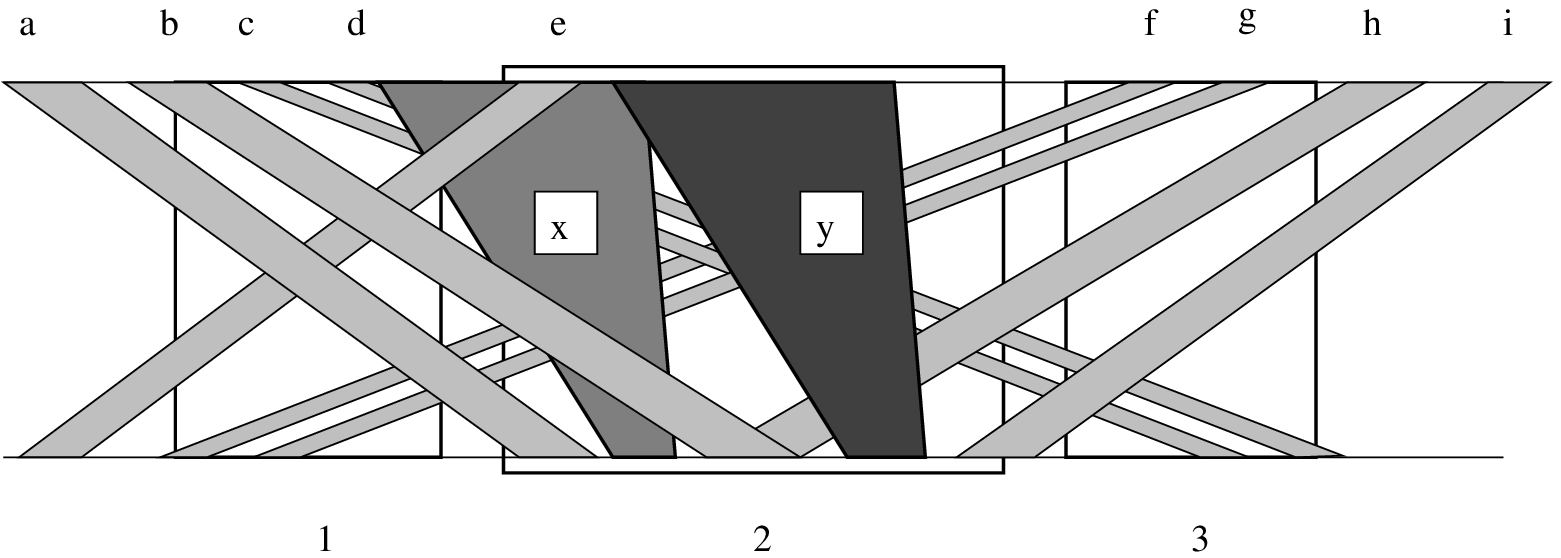}{A generic trapezoid representation of Figure
\ref{fig:uglymess}.  $T_2$ is shown lifted slightly off the baselines
for clarity.}{\hsize}{fig:npput} \fig{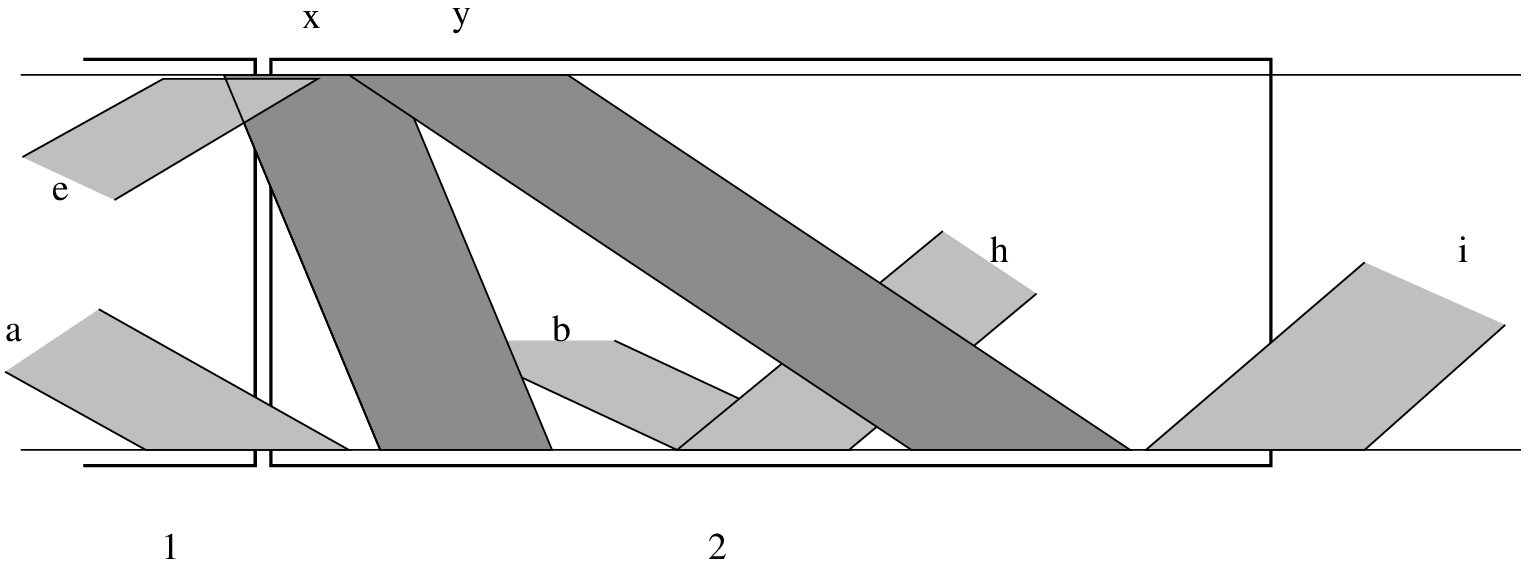}{Detail view of a
parallelogram representation of Figure
\ref{fig:uglymess}.}{\hsize}{fig:npput-p} \fig{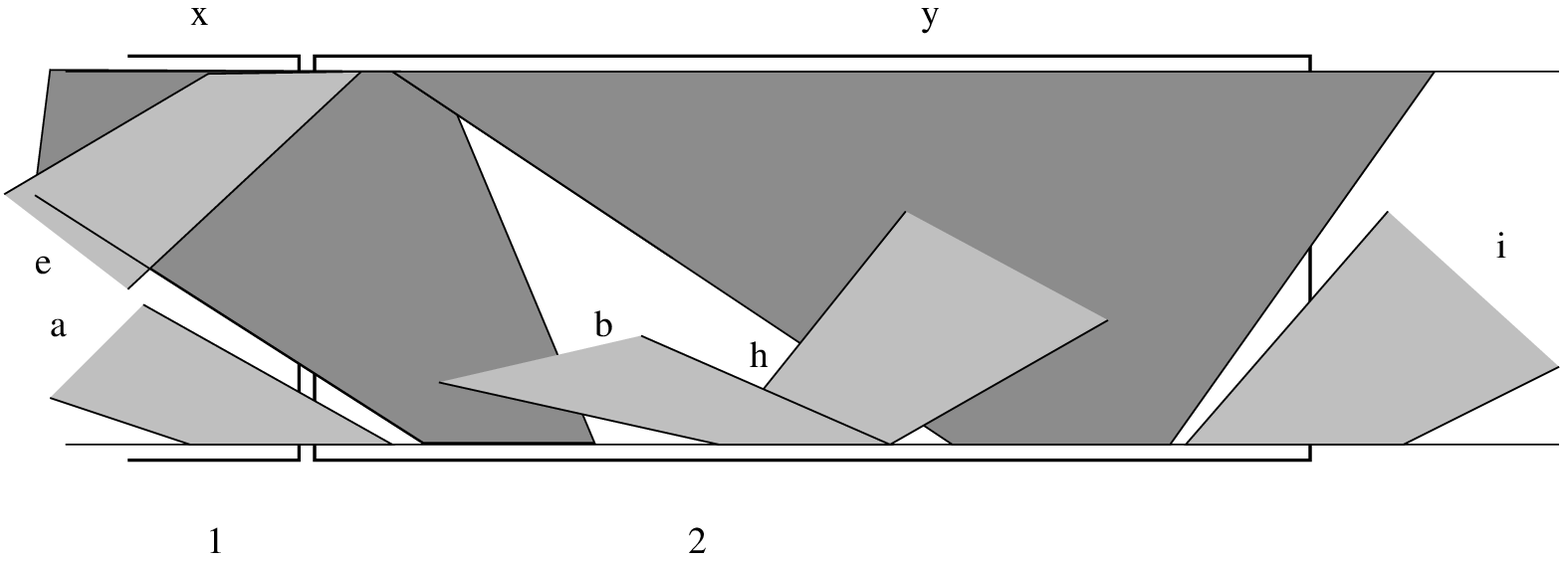}{Detail
view of a unit trapezoid representation of Figure
\ref{fig:uglymess}.}{\hsize}{fig:npput-u}

To show that no representation can be a proper parallelogram
representation, we must first establish the locations of the jaw
structures, so as to have a framework in which to place the offending
parallelograms.  First, we observe that the restrictions to
$A=\{1,2,3,a,c,d,f,g,i\}$,$B=\{1,2,3,b,c,d,f,g,i\}$ and
$C=\{1,2,3,a,c,d,f,g,h\}$ are all copies of the jaw order.

Since $A$ is a jaw, Lemma \ref{lem:jaw} applies to it; without loss of
generality, assume that the jaw formation appears on the lower
baseline.  Note that $B$ and $C$ each differ from $A$ in exactly one
element, one of the teeth of the jaw; therefore, the jaw structures
from the restrictions to $B$ and $C$ must appear on the same baseline
as the jaw structure from $A$.  From applying the jaw lemma to
restriction $B$, we get that $R(b)<L(2)\leq R(2)$; from applying the
jaw lemma to restriction $C$, we get that $R(2)<L(h)$.  Since $b\|h$, we
may conclude that $l(h) \leq r(b)$.  

What this gets us is that the space in between the teeth of the jaw in
restriction $B$ (the ``$a-h$'' jaw) is disjoint from the space between
the teeth of the jaw in restriction $C$ (the ``$b-i$'' jaw).  Since
$x$ is between $a$ and $h$, and $y$ is between $b$ and $i$, this means
that $\Ibot{x}$ and $\Ibot{y}$ are disjoint.  In particular, $\Ibot{x}
\prec \Ibot{y}$, so that $L(y) \leq R(x)$ in order for $x$ to be
incomparable to $y$.

The goal is to show that in any parallelogram representation, $T_y$ is
properly contained in $T_2$.  Since $x$ is incomparable to both $y$
and 1, and $\Ibot{x}$, $\Ibot{1}$ and $\Ibot{y}$ are all pairwise
disjoint, $\Itop{x}$ must intersect both $\Itop{1}$ and $\Itop{y}$.
This will make sure that $\Itop{1}$ and $\Itop{y}$ stay relatively
close together.  This alone, however, does not ensure that $\Itop{y}$
is contained in $\Itop{2}$.  Fortunately, there is one more element in
the order that we have not yet considered, namely, $e$.  Since $e$ is
below $f$, it follows that $\Ibot{e} \prec \Ibot{f} \prec \Ibot{2}$;
therefore $R(e) \geq L(2)$.  Since $y \succ e$, $L(2) < L(y)$.

Nothing established so far depends on the representation being a
parallelogram representation, and so every property established is a
property of every trapezoid representation.  Now assume that the
representation is a parallelogram representation, so that the top and
bottom lengths of every trapezoid are the same.  We already know that
$\Ibot{x}$ and $\Ibot{y}$ are contained in $\Ibot{2}$, and that
$\Ibot{x}$ and $\Ibot{y}$ are disjoint.  It follows from this that
$b(x)+b(y) < b(2)$.  Since we are assuming that the representation is
a parallelogram representation, it is also true that $t(x)+t(y) <
t(2)$.

If we can establish that $R(y)<R(2)$, then $T_y \subset T_2$, and 
any parallelogram representation must be improper.

However, 

$$\begin{aligned}
R(y)&=L(y)+t(y)\\&\leq R(x)+t(y)\\
&=L(x)+t(x)+t(y)\\&<L(2)+t(x)+t(y)\\
&<L(2)+t(2)\\&=R(2)
\end{aligned}
$$

\end{proof}

\section*{Classification}
\fig{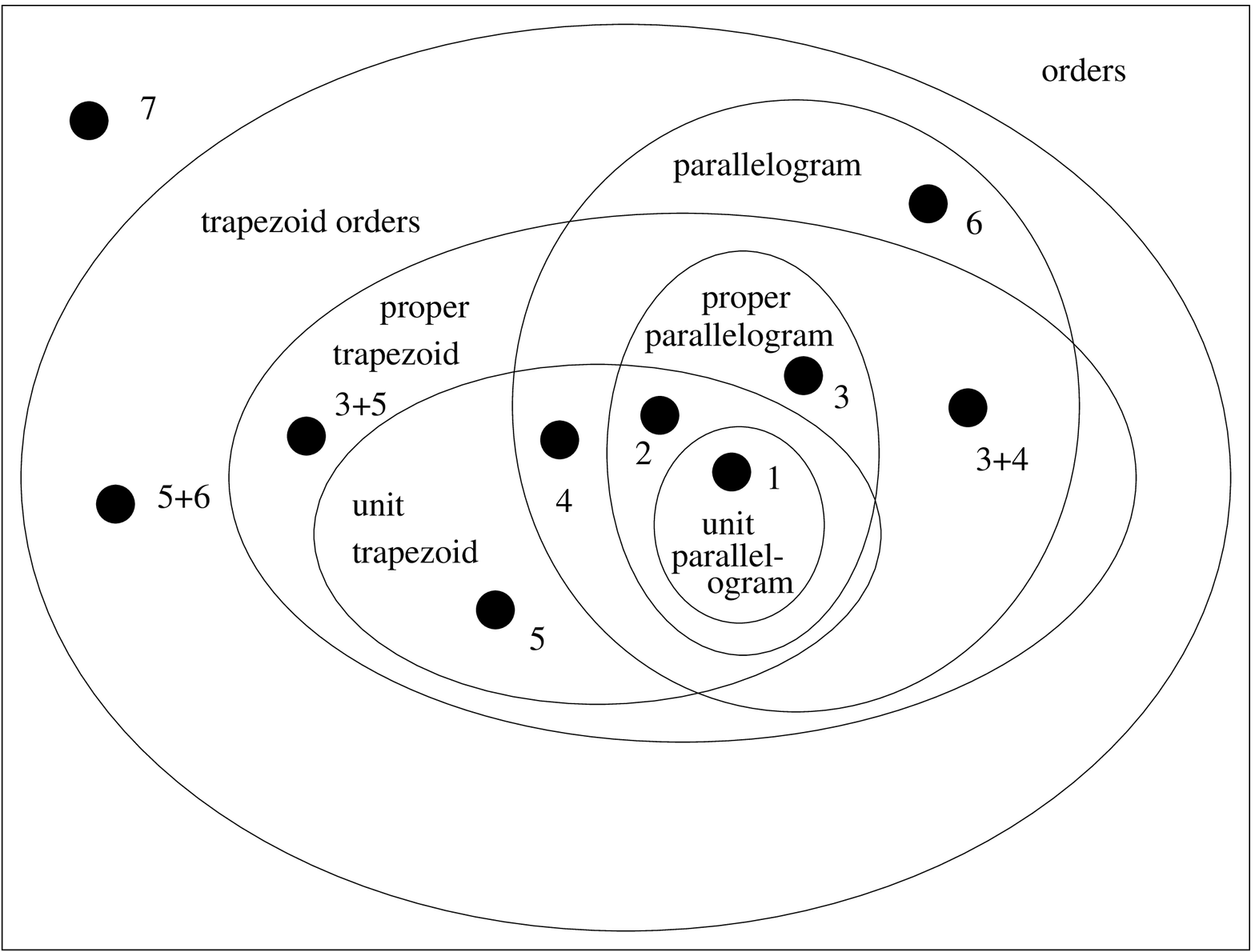}{}{\hsize}{fig:diagram}

The diagram in Figure \ref{fig:diagram} indicates the hierarchy of
different types of trapezoid orders.

In producing the examples for Figure \ref{fig:diagram}, it will be
helpful to recall what we mean by putting two orders together in
series.  If $P$ and $Q$ are ordered sets, then the series sum of $P$
and $Q$ will mean the ordered set whose elements are the disjoint
union of the elements from $P$ and $Q$, and whose comparability
relation is the union of the relations from $P$ and $Q$, together with
the relation that any element in $P$ is less than any element in $Q$.
For our purposes, with trapezoid orders, it is helpful to note that
two orders can be put in series this way by taking trapezoid
representations for each, observing that since both are finite they
are both bounded, and placing the trapezoid representation
of $Q$ past the end of the trapezoid representation of $P$.

The geometric properties that are considered in this paper are
hereditary properties; that is, if an order has one of these
properties, then so does every restriction of that order to any subset
of its elements.  This is easily seen by observing that if a drawing
of the original exists with a particular set of properties, then any
restriction can be shown to have those same properties by taking the
original drawing and removing all trapezoids for elements not in the
restriction.  All trapezoids in the resulting diagram still have the
same properties (and maybe some additional ones), and so the
restriction as a whole has the same properties (and maybe some
additional ones).  As a consequence, if $P$ does not have a particular
property, then neither does the series sum of $P$ and any other order.

Conversely, if two orders both have some property, say having a proper
trapezoid representation, then so does the order formed by taking the
series sum of the two orders; this is easily seen by taking some
representation for each of the two orders having the property in
question, and putting them on two parallel lines as described above.

From this we conclude that the order formed by two orders in series is
contained in the smallest class of orders containing the union of the
smallest class of orders containing each of the two original orders,
but no smaller class.

In Figure \ref{fig:diagram}, the various examples referred to are:

1 - $\underbar 2+\underbar 2$

2 - See Theorem \ref{thm:pput} and Figure \ref{fig:pput}

3 - See Theorem 5 and Figure 13 in \citeBMR

4 - See Theorem \ref{thm:uglymess} and Figure  \ref{fig:uglymess}

5 - See Theorem \ref{thm:nonparallel} and Figure \ref{fig:nonparallel}

6 - See Theorem 2 and Figure 8 in \citeBMR

7 - Any order of interval dimension at least 3; \cite{bogart:intdimbound}
shows that there exist orders of arbitrarily large interval dimension.

\section*{Interval Orders}
Since trapezoid orders are a generalization of interval orders, it is
a natural question to ask where interval orders fit in the hierarchy
described above.

First, we prove a more general lemma on trapezoid orders formed by the
intersection of an interval order and a semiorder.

\begin{lemma}\label{lem:intervalunit}
If $P$ is a trapezoid order formed by the intersection of an interval
order and a semiorder, then $P$ is a unit trapezoid order.
\end{lemma}

\begin{proof}To prove that $P$ is a unit trapezoid order, it is
clearly sufficient to prove that $P$ has a constant area trapezoid
representation, as we can then scale that back down to a unit area
representation.

Let $P = X \cap S$, where $X$ is any interval order and $S$ is a
semiorder on the same set.  Take an interval representation of $X$ and
scale it so that the longest interval in it has length 1.  Since $S$
is a semiorder, it has an interval representation using only unit
length intervals.  Take such a representation and scale it so that no
two distinct endpoints are closer than 2 units apart.  The resulting
intervals all have constant length $k$, because the original interval
representation was a unit-length representation.  Now move the right
endpoint of each interval in $S$ to the right so that the sum of its
length and the length of the corresponding interval in $X$ is $k+1$.

Since the distance between any two endpoints was at least 2, and no
endpoint has been moved more than 1, the relative positions of all
endpoints have been preserved.  Therefore, the new (stretched)
intervals are still a representation of $S$.

The trapezoid representation using the intervals from $X$ and the
stretched intervals from $S$ is composed of trapezoids, the sum of
whose bases is $k+1$, and therefore all these trapezoids have the same area.
\end{proof}

\begin{lemma}\label{lem:intervalproper} If $P$ is a trapezoid order
formed by the intersection of an interval order and a linear order,
then $P$ is a proper parallelogram order.
\end{lemma}

\begin{proof}Let $P=X \cap L$, where $X$ is an interval order and
$L=(x_1<x_2<\ldots<x_n)$ is a linear order.  Take any interval
representation of $X$, and let $l(x)$ stand for the length of the
interval representing $x \in X$.  The following is an interval
representation of $L$:
$[0,l(x_1)],[l(x_1)+1,l(x_1)+l(x_2)+1],
[l(x_1)+l(x_2)+2,l(x_1)+l(x_2)+l(x_3)+2],
\ldots,\left [\left ( \sum_{i=1}^{n-1} l(x_i) \right ) + (n-1),
\left ( \sum_{i=1}^{n} l(x_i) \right ) + (n-1) \right ]$

Each interval has length equal to the length of the corresponding
interval in $X$, and represents a linear order because each interval
is 1 unit away from (and hence does not intersect) any other interval.

The trapezoid representation formed by these two sets of intervals is
a parallelogram representation, because corresponding intervals on the
two lines have the same length.  The trapezoid representation is also
a proper representation, because the bases on one line (the baseline
from the linear order) are all disjoint.  Thus, we have produced a
proper parallelogram representation, so the order must be a proper
parallelogram order.
\end{proof}

\begin{thm}\label{thm:intervalpput}
Any interval order is both a unit trapezoid order and a proper
parallelogram order.
\end{thm}

\begin{proof}Let $I$ be an interval order and $L$ any linear extension
of $I$.  Then $I \cap L = I$ trivially.  Additionally, $L$ may be
viewed as a semiorder, because the class of linear orders is properly
contained in the class of semiorders.

By Lemma \ref{lem:intervalunit}, $I$ must be a unit trapezoid order,
as it can be expressed as the intersection of an interval order and a
semiorder.  By Lemma \ref{lem:intervalproper}, $I$ must be a proper
parallelogram order, as it can be expressed as the intersection of an
interval order and a linear order.
\end{proof}

\begin{lemma}\label{lem:unitpar}
An order is a unit parallelogram order if and only if it has semiorder
dimension at most 2.
\end{lemma}

\begin{proof}Suppose $P$ is an order with semiorder dimension at most
2.  Then there exist two semiorders $S$ and $T$ (possibly identical)
whose intersection is $P$.  Since $S$ and $T$ are semiorders, each has
a unit interval representation.  \cite{scottsuppes} Place the unit
interval representations on parallel lines and construct trapezoids
from them as described in \citeBMR.  Since all
intervals have the same length, namely 1, the top and bottom intervals
for each trapezoid have the same length; i.e. each trapezoid is really
a parallelogram.  Since all intervals have the same length, then the
sum of the top and bottom intervals of each trapezoid is the same;
i.e. each trapezoid has the same area.  Thus $P$ is a unit
parallelogram order.

Suppose now that $P$ is a unit parallelogram order.  Take a unit
parallelogram representation and let $S$ and $T$ be the two interval
orders induced by intersection with the baselines.  Since the
representation is a unit representation, the sum of the bases of each
trapezoid is a constant.  Since the representation is a parallelogram
representation, each base is exactly $\frac 1 2$ of the sum of the
bases, or $\frac 1 2$ of a constant.  Thus, all bases have the same
length.  Since $S$ and $T$ are defined by taking the top and bottom
bases of the trapezoids, all of the intervals in $S$ and $T$ have the
same length; i.e. $S$ and $T$ are unit interval orders.  However, the
property of being a unit interval order is completely equivalent to
the property of being a semiorder, \cite{scottsuppes} and so $P$ is
the intersection of two semiorders, and hence has semiorder dimension
at most 2.
\end{proof}

\begin{lemma}
There exist interval orders with arbitrarily large semiorder
dimension.
\end{lemma}

\begin{proof}
See \cite{bogart:intdimbound} and \cite{trotter:complexity}.  As far
as I know, the result was first stated this way in
\cite{fishburn:semiorders}.
\end{proof}

\fig{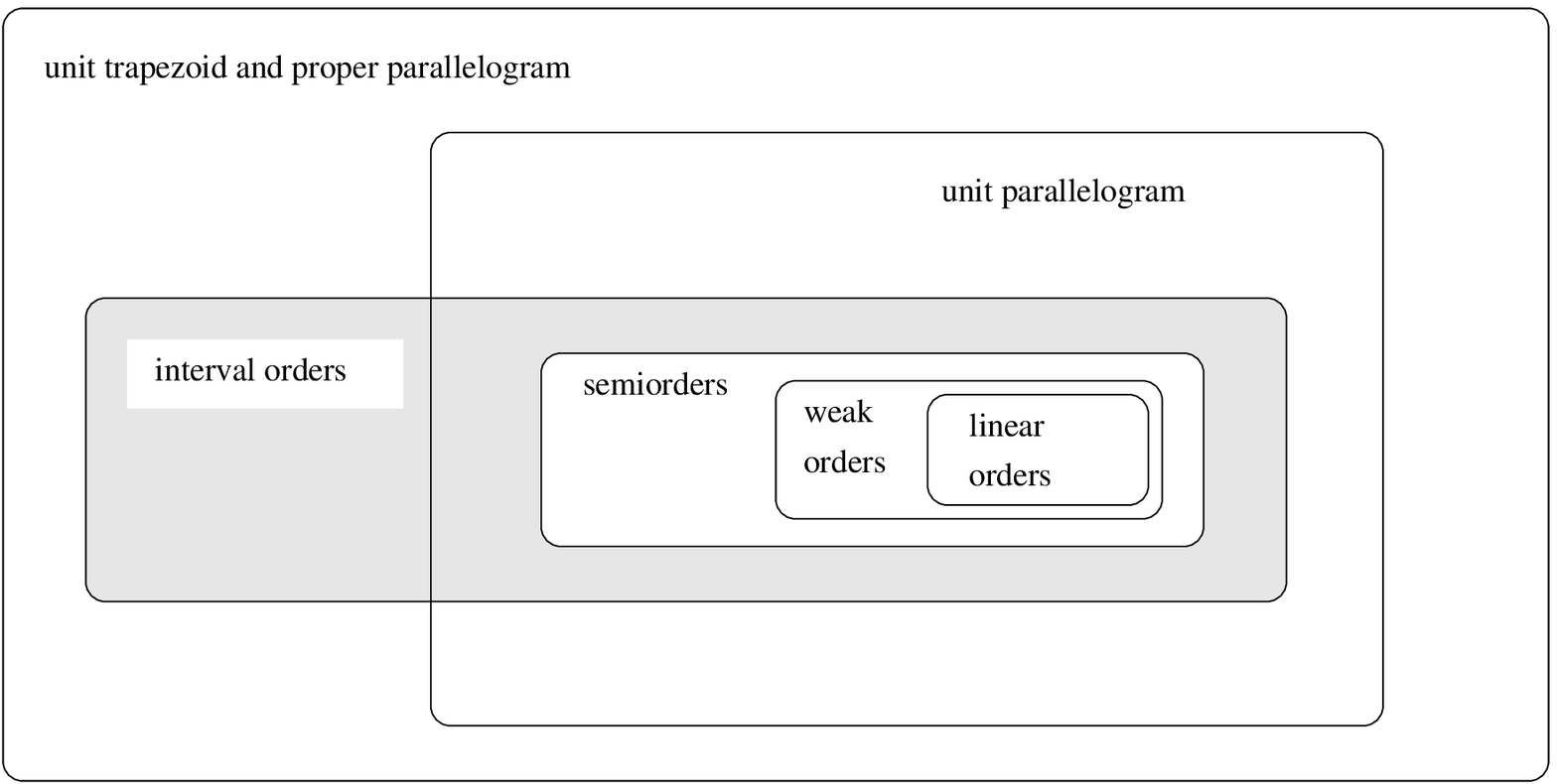}{The positioning of interval orders in the hierarchy
from Figure \ref{fig:diagram}.  These are not shown in the larger
hierarchy in order to keep the picture clear.}
{.75\hsize}{fig:intervaldiagram}

\cite{fishburn85} proves that an order is an interval order if and
only if it has no restriction to a $\underbar 2 + \underbar 2$.  Thus,
our previous example of a unit parallelogram order, $\underbar 2 +
\underbar 2$, shows that there exist unit parallelogram orders which
are not interval orders.

\section*{Application to graphs}
As in \citeBMR, we will now apply the results to trapezoid graphs.
Each class of trapezoid orders gives rise to a corresponding class of
trapezoid graphs, the intersection graphs of the trapezoid
representations of those orders.  Ideally, each property defining a
subclass of trapezoid graphs would be a comparability invariant, so
that each graph arising from an order with that property would also
have that property.  Unfortunately, there are only a few results on
compar ability invariants that apply here.

\cite{gallai67:_trans} establishes the fundamental result in this area,
namely, that one only need consider the autonomous sets, and the effect
of reversing an autonomous set in the order to determine if a property
is a comparability invariant.  

We will begin with the properties that are known to be comparability
invariants.  Habib, Kelly and M\"ohring show that interval dimension
is a comparability invariant in \cite{habib91:_inter}.  This
establishes that trapezoid graphs and interval graphs are precisely
the cocomparability graphs of trapezoid orders and interval orders,
respectively.  In a recent paper, \cite{felsner:_semi} showed that
semi-order dimension two is a comparability invariant.  Thus, the 
property of being a unit parallelogram order is a comparability invariant,
and so the class of unit parallelogram graphs is precisely the class
of cocomparability graphs of unit parallelogram orders.

No further general results are known.  To take care of the remaining
cases, we will fall back on the technique used in \citeBMR.
Fortunately, it is simple to check that the remaining properties we
consider in this paper (proper, unit, parallelogram, and proper
parallelogram) are all comparability invariants of the specific orders
used here.  Each order used as an example has only a few non-trivial
autonomous sets, and they are all short, obvious chains.  Each of
these autonomous chains will actually result in an isomorphic order
when reversed, and so the properties will be the same.  The arguments
used to show this are not included here, as they are repetitions of
the arguments used in Proposition 3 of \citeBMR, only longer due
to the loss of some of the symmetry used there.

\bibliographystyle{apalike}
\bibliography{ryan}
\end{document}